\documentclass[12pt,a4paper]{article}
\usepackage[english]{babel}
\usepackage{latexsym,bm}
\begin{document}
\newtheorem{axiom}{Axiom}
\newtheorem{definition}{Definition}
\newtheorem{theorem}{Theorem}
\pagestyle{plain}
\title{The Foundations of Mathematics\\in the Physical Reality}
\author{Doeko H. Homan}
\date{May 30, 2022}
\maketitle
\pagenumbering{arabic}
It is well-known that the concept {\em set} can be used as the foundations for
mathematics. And the Peano axioms for the set of all natural numbers should be
`considered as the fountainhead of all mathematical knowledge' (Halmos [1974]
page 47). However, natural numbers should be {\em defined}, thus `what is a
natural number', not `what is the set of all natural numbers'.\par
The basic properties of a set are the members of that set, and the sets that
set is a member of. Thus there is no set without members. And set theory has to
be included objects existing in physical reality, called `individuals'.
Example: A single shoe is not a singleton set but an individual, a pair of
shoes (often a left one and a right one) is a set of individuals.\par
In this article we present an axiomatic definition of sets with individuals.
Natural numbers and ordinals are defined. Limit ordinals are `first numbers',
that is a first number of the Peano axioms. For every natural number $m$ we
define `first $\omega^m$-numbers'. Every ordinal is an ordinal with a first
$\omega^0$-number, every `limit ordinal' is an ordinal with a first
$\omega^m$-number. Ordinals with a first number satisfy the Peano axioms.
\section{What is a set?}
At an early age you develop the idea of `what is a set'. You belong to a family
or a clan, you belong to the inhabitants of a village or a particular region.
Experience shows there are objects that constitute a set. A constituent of a
set is called `a member' of that set. A set can be a member of a set.\par
Experience shows an object `is equal to' or `is not equal to' another object.
The relations `is equal to' and `is not equal to' are defined in set theory.
For individuals the relations `is equal to' and `is not equal to' are given in
physical reality but are not known within set theory. Therefore we define: The
only member of an individual is the individual oneself.
\section{Very basic set theory}
Sets are denoted by $s$, $t$, $u$, $\ldots$, by capital letters or by lowercase
Greek letters. Symbol $\in$ is the `membership relation symbol', to read as `is
a member of' or as `belongs to'. For every $s$ and every $t$ applies $s\in t$
or the negation of $s\in t$, denoted by $s\notin t$. Thus for every $s$ applies
$s\in s$ or $s\notin s$.\par
Formulas are built up from formulas like $s\in t$ by means of $\wedge$ (and),
$\vee$ (or), $\rightarrow$ (material implication), $\leftrightarrow$
(biimplication), $\neg$ (negation) and by means of $\forall s$ (for every set
$s$) and $\exists s$ (exists set $s$). The parentheses $($ and $)$ preclude
ambiguity and improve legibility. Then $s\notin t$ is short for $\neg(s\in t$),
$(s\in t\vee s\notin t)$ is a tautology and $\exists s\exists t(s\in t\wedge s
\notin t)$ is a contradiction.
\begin{definition}
Set $s$ is an `individual' if $s\in s$.
\hfill\boldmath$\Box$\end{definition}
If $s$ is an individual then $s$ `is a member of oneself'. Sets being a member
of themselves seem problematic. Section 4.6.5 in (Mendelson [2015]) states:\par
`The words "individual" and "atom" are sometimes used as synonyms for\par
"urelement". (...) Thus, urelements have no members.'\\
Thus an atom or an urelement is certainly not what we call an individual. The
membership relation for individuals seems awkward. But there are no logical
obstacles for `being a member of oneself'. To make a distinction between
individuals and other sets, in the phrase `an individual or a set', or `a set
of individuals', or `a transitive set', that sets are not individuals.\\
\mbox{ }\par
Then $\forall s(s\in s\vee s\notin s)$ thus an individual is a member of $s$,
or a set is not a member of $s$. Therefore $\forall s\exists u((u\in s\wedge u
\in u)\vee(u\notin s\wedge u\notin u))$.\\
The same result also follows from Russell's paradox.\\
\mbox{ }\\
Russell's paradox (Jech [2002] page 4).\\
Assume $s$ is a set whose members are all those (and only those) sets that are
not members of themselves\par
$\exists s\forall u(u\in s\leftrightarrow u\notin u)$.\\
The formula applies to every $u$ thus we can choose $s$. Then the formula
reads\par
$\exists s(s\in s\leftrightarrow s\notin s)$ thus $\exists s(s\in s\wedge s
\notin s)$, a clear contradiction. Therefore\par
$\forall s\exists u((u\in s\wedge u\in u)\vee(u\notin s\wedge u\notin u))$.\\
Thus Russell's paradox shows: For every $s$, an individual is a member of $s$,
or a set is not a member of $s$. Logic allows for individuals.\\
\mbox{ }\\
We try to prove by the same kind of reasoning $\exists v(v\notin s)$.\\
Assume $v$ is a set whose members are all those (and only those) members of $s$
that are not members of themselves (thus are not individuals)\par
$\exists v\forall u(u\in v\leftrightarrow(u\in s\wedge u\notin u))$.\\
The formula applies to every $u$ thus we can choose $v$\par
$\exists v(v\in v\leftrightarrow(v\in s\wedge v\notin v))$.\\
Then $\exists v(v\notin v\wedge v\notin s)$ thus $v$ is not an individual and
not a member of $s$. Thus (Jech [2002] page 4): `The set of all sets does not
exist.'\\
\mbox{ }\\
Again we try the same kind of reasoning. Assume $v$ is a set whose members are
all those (and only those) members of $s$ that are individuals\par
$\exists v\forall u(u\in v\leftrightarrow(u\in s\wedge u\in u))$.\\
The formula applies to every $u$ thus we can choose $v$\par
$\exists v(v\in v\leftrightarrow(v\in s\wedge v\in v))$.\\
Then $\exists v(v\notin v\vee(v\in v\wedge v\in s))$. Thus $v$ is a set of
individuals, or $v$ is an individual and a member of $s$. Again, individuals
are on par with other sets. And `the set of all individuals' is not excluded.
\section{The axioms of set theory}
\begin{definition} $s$ `is equal to' $t$, denoted by $s=t$ if $\forall u(u\in s
\leftrightarrow u\in t)$.\par
$s\neq t$ is short for $\neg(s=t)$.
\hfill\boldmath$\Box$\end{definition}
\begin{axiom}(equality)
$\forall s\forall t(s=t\rightarrow\forall u(s\in u\leftrightarrow t\in u))$.
\hfill\boldmath$\Box$\end{axiom}
If $s=t$ then $s$ and $t$ have the same basic properties: $s$ and $t$ have the
same members, and $s$ and $t$ are members of the same sets.
\begin{axiom}(individuals)
$\forall s(s\in s\rightarrow\forall u(u\in s\rightarrow u=s))$.
\hfill\boldmath$\Box$\end{axiom}
Thus the only member of an individual is the individual oneself. Then\par
$\forall s\forall t(s\in s\wedge t\notin t\rightarrow s\neq t\wedge t\notin
s)$.
\begin{axiom}(pairs)
$\forall s\forall t\exists v\forall u(u\in v\leftrightarrow(u=s\vee u=t))$.
\hfill\boldmath$\Box$\end{axiom}
Then $v$ is denoted by $\{s,t\}$. The `singleton' $\{s\}$ of $s$ is short for
$\{s,s\}$.\par
$s\notin s\rightarrow\{s\}\neq s$ otherwise $s\in s$. Therefore $s\notin s
\rightarrow\{s\}\notin\{s\}$.\par
$s\in s\rightarrow\{s\}=s$, the singleton of an individual is the individual
oneself.\par
$s\neq t\rightarrow\{s,t\}\notin\{s,t\}$. Thus a pair of different sets is not
an individual.\\
\mbox{ }\\
The `empty set' is usually postulated by $\exists s\forall u (u\notin s)$. The
formula applies to every $u$ thus we can choose $s$. Therefore $\exists s(s
\notin s)$ thus the empty set, denoted by $\emptyset$, is not an individual.
The axiom of regularity is usually\par
$\forall s(s\neq\emptyset\rightarrow\exists v(v\in s\wedge\forall u(u\in v
\rightarrow u\notin s)))$.\\
Thus every individual contradicts the axiom of regularity. Therefore we have to
abandon the empty set and reformulate the axiom of regularity.
\begin{axiom} $\forall s\exists v(v\in s\wedge v\in v\wedge\forall u(u\in s
\rightarrow u\in u)\vee$\\
(regularity) $v\in s\wedge v\notin v\wedge\forall u(u\in v\wedge u\notin u
\rightarrow u\notin s))$.
\hfill\boldmath$\Box$\end{axiom}
Then to every set belongs an individual or a set, and regularity does not apply
to $s$ if $s$ is an individual or a set of individuals.
\begin{axiom}(union)
$\forall s\exists v\forall u(u\in v\leftrightarrow\exists w(w\in s\wedge u\in
w))$.
\hfill\boldmath$\Box$\end{axiom}
The union $v$ of $s$ is denoted by $\bigcup s$, and $\forall s\exists u(u\in
\bigcup s)$. The `union $s\cup t$' of $s$ and $t$ is defined by $s\cup t=
\bigcup\{s,t\}$. Thus $s\in s\rightarrow\bigcup s=s\wedge s\cup\{s\}=s$.
\begin{axiom}(power)
$\forall s\exists v\forall u(u\in v\leftrightarrow\forall w(w\in u\rightarrow
w\in s))$.
\hfill\boldmath$\Box$\end{axiom}
The power $v$ of $s$ is denoted by ${\cal{P}}s$, and $\forall s(s\in
{\cal{P}}s)$.\\
We reformulate the axiom `separation'.
\begin{axiom}(separation)\par
If $\Phi(u)$ denotes a well formed formula $\Phi$ with $u$ free in $\Phi$
then\par
$\forall s(\exists u(u\in s\wedge\Phi(u))\rightarrow\exists v\forall u(u\in v
\leftrightarrow u\in s\wedge\Phi(u)))$.
\hfill\boldmath$\Box$\end{axiom}
In section 7 we postulate an `axiom of infinity'.
\section{Transitive sets}
\begin{definition}
$s$ is `included in' $t$, denoted by $s\subseteq t$ if $\forall u(u\in s
\rightarrow u\in t)$.\par
$s$ `is transitive' if $\forall u(u\in s\rightarrow u\subseteq s)$, or
equivalently $\bigcup s\subseteq s$.
\hfill\boldmath$\Box$\end{definition}
Every individual and every set of individuals is transitive. If $s$ is
transitive then $\bigcup s$ and $s\cup\{s\}$ are transitive. If $v\notin v$
then $\{v\}$ is not transitive.
\begin{theorem} If $s$ is transitive and $\exists u(u\in s\wedge u\notin u)$
then a set of individuals\par
is a member of $s$: $\exists v(v\in s\wedge v\notin v\wedge\forall u(u\in v
\rightarrow u\in u))$.
\end{theorem}
{\it Proof.} Apply regularity to $s$ to find $v\in s\wedge v\notin v$ such
that\par
$\forall u(u\in v\wedge u\notin u\rightarrow u\notin s)$. Then $\forall u(u\in
v\wedge u\in s\rightarrow u\in u)$.\par
Set $s$ is transitive and $v\in s$ thus $\forall u(u\in v\rightarrow u\in
s)$.\par
Therefore $\exists v(v\in s\wedge v\notin v\wedge\forall u(u\in v\rightarrow u
\in u))$.
\begin{flushright}\boldmath$\Box$\end{flushright}
If $u$ and $v$ are the only individuals belonging to transitive set $s$ then
$\{u,v\}$ is the only possible set of individuals belonging to $s$. Therefore
$\{u,v\}\in s\cup\{s\}$.\par
In the remainder we fix the pair of different individuals $o$ and $e$ and
define $0=\{o,e\}$. Then $0\notin0$, $\bigcup0=0$ and $\forall u(u\notin u
\rightarrow u\notin0)$.
\begin{theorem}
If $T$ is a transitive set with transitive members and\par
$0\in T\wedge\forall u(u\in T\wedge u\in u\rightarrow u\in0)$ then\par
$\forall s\forall t(s\in T\wedge s\notin s\wedge t\in T\wedge t\notin t
\rightarrow(s\in t\vee s=t\vee t\in s))$.
\end{theorem}
{\it Proof.} $T$ is a transitive set with transitive members and $0\in T$
therefore\par
$\forall u(u\in T\rightarrow((u\notin0\leftrightarrow u\notin u)\wedge(u\notin
u\rightarrow0\in u\cup\{u\})))$. Thus\par
$0\in T\wedge0\notin0\wedge\forall u(u\in T\wedge u\notin u\rightarrow(0\in u
\vee0=u\vee u\in0))$.\\
Assume exists $v$ such that\par
$v\in T\wedge v\notin v\wedge\exists w(w\in T\wedge w\notin w\wedge(v\notin w
\wedge v\neq w\wedge w\notin v))$.\par
Then $v\neq0$ and $v\neq w$.\\
Use separation to define $V$ by\\
$\forall u(u\in V\leftrightarrow u\in T\wedge u\notin u\wedge\exists x(x\in T
\wedge x\notin x\wedge(u\notin x\wedge u\neq x\wedge x\notin u)))$.\par
Then $v\in V\wedge v\notin v$, $0\notin V$ and $\forall u(u\in0\rightarrow u
\notin V)$.\\
Apply regularity to $V$ to find $t\in V$ such that $t\in T\wedge t\notin t
\wedge$\par
$\exists x(x\in T\wedge x\notin x\wedge(x\notin t\wedge x\neq t\wedge t\notin
x))\wedge\forall u(u\in t\wedge u\notin u\rightarrow u\notin V)$.\par
$0\notin V$ thus $t\neq0$. And $t\in T\wedge t\notin t$ therefore $0\in t$.\par
Set $t\in V$ thus $\exists x(x\in T\wedge x\notin x\wedge(x\notin t\wedge x\neq
t\wedge t\notin x))$.\\
Use separation to define $W$ by\par
$\forall u(u\in W\leftrightarrow u\in T\wedge u\notin u\wedge(u\notin t\wedge u
\neq t\wedge t\notin u))$.\par
Then $x\in W\wedge x\notin x$ and $0\notin W$.\\
Apply regularity to $W$ to find $y\in W$ such that\par
$y\in T\wedge y\notin y\wedge(y\notin t\wedge y\neq t\wedge t\notin y)\wedge
\forall u(u\in y\wedge u\notin u\rightarrow u\notin W)$.\par
Then $y\neq0\wedge y\in T\wedge y\notin y$ thus $y$ is transitive and $0\in
y$.\par
Set $T$ is transitive thus $\forall u(u\in y\wedge y\in T\rightarrow u\in T)$
therefore\\
$\forall u(u\in y\wedge u\notin u\rightarrow u\notin W)\rightarrow\forall u(u
\in y\wedge u\notin u\rightarrow(u\in t\vee u=t\vee t\in u))$.\par
Then $(t\in u\vee t=u)\wedge u\in y\rightarrow t\in y$ contradictory to $t
\notin y$.\par
Thus $\forall u(u\in y\wedge u\notin u\rightarrow u\notin W)\rightarrow
\forall u(u\in y\wedge u\notin u\rightarrow u\in t)$.\par
Sets $y$ and $t$ are transitive, $0\in y$ and $0\in t$ thus $y\subseteq t$.\par
Then $y\subseteq t\wedge y\neq t\rightarrow\exists s(s\in t\wedge s\notin s
\wedge s\notin y)$ otherwise $y=t$.\par
Set $s\neq0$, $s\in t$ and $\forall u(u\in t\wedge u\notin u\rightarrow u\notin
V)$ therefore\par
$\forall u(u\in T\wedge u\notin u\rightarrow(u\in s\vee u=s\vee s\in u))$.\par
Set $y\in T$ and $s\notin y$ thus $y\in s\vee y=s$.\par
But $(y\in s\vee y=s)\wedge s\in t\rightarrow y\in t$ contradictory to $y\notin
t$.\\
Therefore the assumption `exists $v$' is not correct thus\par
$\forall s\forall t(s\in T\wedge s\notin s\wedge t\in T\wedge t\notin t
\rightarrow(s\in t\vee s=t\vee t\in s))$.
\begin{flushright}\boldmath$\Box$\end{flushright}
If both $s$ and $t$ are transitive sets with transitive members then $(s\cup t)
\cup\{s,t\}$ is a transitive set with transitive members and $s$ and $t$ belong
to $(s\cup t)\cup\{s,t\}$. Then the `law of trichotomy' reads\par
$s\notin s\wedge t\notin t\wedge\forall u(u\in s\cup t\wedge u\in u\rightarrow
u\in0)\rightarrow s\in t\vee s=t\vee t\in s$.\\
Then $\bigcup s=s$ or $\bigcup s\neq s\wedge s=\bigcup s\cup\{\bigcup s\}$.
\section{Natural numbers}
\begin{definition}
$s$ is a `natural number' if $s$ is a transitive set with transitive\par
members and $s\notin0\wedge\forall u(u\in s\cup\{s\}\rightarrow u\in0\cup\{0\}
\vee\bigcup u\in u)$.\par
The `successor' of natural number $s$ is $s\cup\{s\}$.
\hfill\boldmath$\Box$\end{definition}
Then $0$ is the first natural number. In fact, any pair of different
individuals can be a first natural number. In this section $0$ is denoted by
$\alpha$. If $s$ is a natural number then $s\notin\alpha\leftrightarrow s\notin
s$. And
$$\bigcup\alpha=\alpha\wedge s\notin\alpha\wedge\forall u(u\in s\cup\{s\}
\rightarrow u\in\alpha\cup\{\alpha\}\vee\bigcup u\in u).$$
Then $\alpha$ is the first natural number. And\\
\mbox{ }\par
$\forall u(u\in s\cup\{s\}\wedge u\in u\rightarrow u\in\alpha)$,\par
$\forall u(u\in s\cup\{s\}\rightarrow u\in\alpha\cup\{\alpha\}\vee u=\bigcup u
\cup\{\bigcup u\})$,\par
$\forall u(u\in s\cup\{s\}\wedge u\notin\alpha\rightarrow\bigcup u\notin
\alpha)$.
\begin{theorem}
The successor $s\cup\{s\}$ of natural number $s$ is a natural number.
\end{theorem}
{\it Proof.} $s\cup\{s\}$ is a transitive set with transitive members.\par
$s$ is a natural number thus $\forall u(u\in s\cup\{s\}\rightarrow u\in\alpha
\cup\{\alpha\}\vee\bigcup u\in u)$.\par
And $s\notin\alpha$ implies $s\cup\{s\}\notin\alpha$ and $s\cup\{s\}\neq
\alpha$. Thus $s\cup\{s\}\notin\alpha\cup\{\alpha\}$.\par
Then $\bigcup(s\cup\{s\})=s\wedge s\in s\cup\{s\}$ thus $\bigcup(s\cup\{s\})\in
s\cup\{s\}$.\par
Therefore $s\cup\{s\}$ is a natural number.
\begin{flushright}\boldmath$\Box$\end{flushright}
\begin{theorem}
If $s$ is a natural number and $\alpha\in s$ then\par
$\forall u(u\notin\alpha\wedge u\in s\rightarrow(u\mbox{ is a natural
number}))$.
\end{theorem}
{\it Proof.} If $u\in s$ then $u$ is a transitive set with transitive
members.\par
$s$ is transitive thus $\forall v(v\in u\cup\{u\}\wedge u\in s\rightarrow v\in
s)$. Then\par
$\bigcup\alpha=\alpha\wedge u\notin\alpha\wedge\forall v(v\in u\cup\{u\}
\rightarrow v\in\alpha\cup\{\alpha\}\vee\bigcup v\in v)$.\par
Therefore $\forall u(u\notin\alpha\wedge u\in s\rightarrow(u\mbox{ is a natural
number}))$.
\begin{flushright}\boldmath$\Box$\end{flushright}
Next theorem is the `principle of mathematical induction'.
\begin{theorem}
If $s$ is a natural number and $\Phi(u)$ is a well formed formula $\Phi$\par
with natural number $u$ free in $\Phi$ then\par
$\Phi(\alpha)\wedge\forall u(u\notin\alpha\wedge u\in s\wedge\Phi(u)\rightarrow
\Phi(u\cup\{u\}))\rightarrow\Phi(s)$.
\end{theorem}
{\it Proof.} If $s$ is first number $\alpha$ then $\Phi(\alpha)$ thus
$\Phi(s)$.\par
If $s\neq\alpha$ then $\alpha\in s$. Use separation to define set $v$ by\par
$\forall u(u\in v\leftrightarrow u\notin\alpha\wedge u\in s\wedge\Phi(u))$.
Then $\alpha\in v$ and $s\notin v$.\par
Use separation to define set $w$ by\par
$\forall u(u\in w\leftrightarrow u\notin\alpha\wedge u\in s\cup\{s\}\wedge
u\notin v)$. Then $s\in w$ and $\alpha\notin w$.\par
Apply regularity to set $w$ to find natural number $z\in w$ such that\par
$z\in s\cup\{s\}\wedge z\notin\alpha\wedge z\notin v\wedge\forall u(u\in z
\wedge u\notin u\rightarrow u\notin w)$.\par
Then $\alpha\in v$ thus $z\neq\alpha$. And $z\notin\alpha$ therefore $\alpha\in
z$. Thus $z$ is not first\par
number $\alpha$ therefore $z=\bigcup z\cup\{\bigcup z\}$. Then $\bigcup z\in z$
and $\bigcup z\notin\alpha$.\par
And $\bigcup z\in z\wedge z\in s\cup\{s\}\rightarrow\bigcup z\in s$ thus
$\bigcup z\in s\wedge\bigcup z\notin w\rightarrow\bigcup z\in v$.\par
Therefore $\Phi(\bigcup z)$ thus $\Phi(z)$. And $z\in s\cup\{s\}$.\par
If $z\in s$ then $z\notin\alpha\wedge\Phi(z)\rightarrow z\in v$ contradictory
to $z\notin v$.\par
Therefore $z=s$ thus $\Phi(s)$.
\begin{flushright}\boldmath$\Box$\end{flushright}
The axioms formulated by G. Peano are satisfied:
\begin{itemize}
\item there is first natural number $\alpha$,
\item every natural number has as successor a natural number,
\item if $s$ is a natural number and $s\cup\{s\}=\alpha$ then $s\in\alpha$
      contradictory to $s\notin\alpha$ thus first number $\alpha$ is not the
      successor of a natural number,
\item if $s$ and $t$ are natural numbers and $s\cup\{s\}=t\cup\{t\}$ then
      $s\in t\cup\{t\}$ and $t\in s\cup\{s\}$ therefore $s=t$ thus if both
      successors of a pair of natural numbers are equal then both natural
      numbers are equal,
\item the principle of mathematical induction.
\end{itemize}
In the remainder natural numbers are denoted by $l$, $m$, $n$, $x$ or $y$.
The first natural number is $0$, and $1$ denotes the successor of $0$ thus $1=0
\cup\{0\}$.
\begin{definition}
$m$ is `less than' $n$, denoted by $m<n$ if $m\in n$,\par
$m$ is `less than or equal to' $n$, denoted by $m\leq n$ if $n\notin m$.
\hfill\boldmath$\Box$\end{definition}
Then $0$ is less than every other natural number. We define addition $m+n$ and
multiplication $m\cdot n$. As usual multiplication precedes addition.
\begin{definition}
$m+0=m\mbox{, }\mbox{ }m+(n\cup\{n\})=(m+n)\cup\{m+n\}$,\par
$m\cdot0=0\mbox{, }\mbox{ }m\cdot(n\cup\{n\})=m\cdot n+m$.
\hfill\boldmath$\Box$\end{definition}
Then $m+1=(m+0)\cup\{m+0\}$ thus $m+1=m\cup\{m\}$ and $0+1=1$. With the
principle of mathematical induction one can prove $m+n$ and $m\cdot n$ are
natural numbers, $1+n=n+1$ and the basic properties of the Peano Arithmetic
(addition and multiplication).
\section{Ordinals}
\begin{definition}
$\beta$ is an `ordinal' if $\beta$ is a transitive set with transitive\par
members and $\beta\notin0\wedge\forall u(u\in\beta\wedge u\in u\rightarrow u
\in0)$.\par
Ordinal $\alpha$ is a `first number' if $\bigcup\alpha=\alpha$.
\hfill\boldmath$\Box$\end{definition}
Thus $0\in\beta\cup\{\beta\}$ and $0$ is a first number. If $\alpha$ is a first
number and $u\in\alpha\wedge u\notin u$ then $u$ and $u\cup\{u\}$ are ordinals
and $u\cup\{u\}\in\alpha$. And $\beta\notin0\leftrightarrow\beta\notin\beta$.
We define addition of ordinals and natural numbers.
\begin{definition}
If $\beta$ is an ordinal then\par
$\beta+0=\beta\mbox{, }\mbox{ }\beta+(n+1)=(\beta+n)\cup\{\beta+n\}$.
\hfill\boldmath$\Box$\end{definition}
Then $\beta\cup\{\beta\}$ and $\beta+1$ are ordinals. With the principle of
mathematical induction one can prove $\forall m\forall n((\beta+m)+n=\beta+
(m+n))$.
\begin{theorem}
If $\beta$ is an ordinal then $\exists\alpha\exists n(\bigcup\alpha=\alpha
\wedge\beta=\alpha+n)$.
\end{theorem}
{\it Proof.} If $\beta$ is a first number then $\beta=\beta+0$ thus $\alpha=
\beta$ and $n=0$.\par
If $\bigcup\beta\neq\beta$ then $\beta=\bigcup\beta+1$. Use separation to
define set $v$ by\par
$\forall u(u\in v\leftrightarrow u\in\beta\wedge u\notin u\wedge\exists n
(\beta=u+n))$. Then $\bigcup\beta\in v$.\par
Apply regularity to $v$ to find ordinal $\alpha\in v\wedge\alpha\notin\alpha$
such that\par
$\alpha\in\beta\wedge\exists n(\beta=\alpha+n)\wedge\forall u(u\in\alpha\wedge
u\notin u\rightarrow u\notin v)$. Then $\bigcup\alpha\in\alpha\cup
\{\alpha\}$.\par
If $\bigcup\alpha\in\alpha$ then $\alpha=\bigcup\alpha+1$ therefore $\beta=
\bigcup\alpha+(1+n)$.\par
Then $\bigcup\alpha\in\alpha\wedge\alpha\in\beta$ thus $\bigcup\alpha\in\beta$
therefore $\bigcup\alpha\in v$.\par
Then $\bigcup\alpha\in\alpha\wedge\bigcup\alpha\in v$ contradictory to $\forall
u(u\in\alpha\wedge u\notin u\rightarrow u\notin v)$.\par
Therefore $\exists\alpha\exists n(\bigcup\alpha=\alpha\wedge\beta=\alpha+n)$.
\begin{flushright}\boldmath$\Box$\end{flushright}
\begin{definition}
Ordinal $\beta$ is an `ordinal with first number $\alpha$' if\par
$\bigcup\alpha=\alpha\wedge\beta\notin\alpha\wedge\forall u(u\in\beta\cup
\{\beta\}\rightarrow u\in\alpha\cup\{\alpha\}\vee\bigcup u\in u)$.\par
The `successor' of ordinal $\beta$ is $\beta\cup\{\beta\}$.
\hfill\boldmath$\Box$\end{definition}
Then $\forall\beta(\beta\notin\alpha\rightarrow\beta\notin\beta)$.
First number $\alpha$ satisfies the definition of an ordinal with first number
$\alpha$. Natural numbers (the `finite' numbers) are ordinals with first number
$0$.\par
The theorems in section 5 about natural numbers apply analogous to ordinals
with first number $\alpha$. Thus the successor of an ordinal with first number
$\alpha$ is an ordinal with first number $\alpha$ and the ordinals with first
number $\alpha$ satisfy the Peano axioms.\par
If $\beta$ is an ordinal with first number $\alpha$ and $\beta\neq\alpha$ then
$\bigcup\beta\in\beta$ thus $\bigcup\beta$ is the `greatest' ordinal with first
number $\alpha$ belonging to $\beta$.
\begin{theorem}
If $\beta$ is an ordinal with first number $\alpha$ then for every $n$ $(\beta+
n)$\par is an ordinal with first number $\alpha$.
\end{theorem}
{\it Proof.} $\beta+0=\beta$ thus $\beta+0$ is an ordinal with first number
$\alpha$.\par
Assume if $x<n$ then $\beta+x$ is an ordinal with first number $\alpha$. Then
the\par successor $(\beta+x)+1$ of ordinal $\beta+x$ is an ordinal with first
number $\alpha$\par
therefore $\beta+(x+1)$ is an ordinal with first number $\alpha$.\par
Apply the principle of mathematical induction to $0\leq x\wedge x<n$ and\par
conclude for every $n$ $(\beta+n)$ is an ordinal with first number $\alpha$.
\begin{flushright}\boldmath$\Box$\end{flushright}
Every ordinal is an ordinal with a first number $\alpha$. We define $\omega^0=
1$. Then\\\mbox{ }\par
$\forall n(\alpha+\omega^0\cdot n$ is an ordinal with first number
$\alpha)$,\par
$\forall n(\alpha+\omega^0\cdot(n+1)=(\alpha+\omega^0\cdot n)+\omega^0)$,\par
$\forall n(\alpha+\omega^0\cdot n\in(\alpha+\omega^0\cdot n)+\omega^0)$.
\section{More first numbers}
\begin{definition}
If $\gamma$ is an ordinal then\par
$F^0_\gamma=\gamma\mbox{, }\mbox{ }\forall u(u\in F^{m+1}_\gamma\leftrightarrow
u\in F^m_\gamma\wedge\bigcup F^m_u=u)$.\par
Ordinal $\alpha$ is a `first $\omega^m$-number' if $\bigcup F^m_\alpha=\alpha$.
\hfill\boldmath$\Box$\end{definition}
Then $\forall m(F^m_0=0\wedge\bigcup F^m_0=0\wedge F^m_1=1\wedge\bigcup F^m_1=0
\wedge(0\in F^m_\gamma\cup\{F^m_\gamma\}))$.\\
If $0<m$ then $F^m_\gamma$ is a set of individuals and first numbers belonging
to $\gamma$.\\
And $\bigcup\bigcup F^m_\gamma\subseteq\bigcup F^m_\gamma$. If $u\in\bigcup
F^m_\gamma$ then $u$ is transitive and $u\cup\{u\}\in\bigcup F^m_\gamma$\par
thus $u\in\bigcup\bigcup F^m_\gamma$. Therefore if $0<m$ then $\bigcup
F^m_\gamma$ is a first number.\\
Then $\forall u(u\in F^m_\gamma\rightarrow\bigcup F^m_\gamma\notin u)$.\\
Apply the law of trichotomy and conclude $0<m\rightarrow F^m_\gamma\subseteq
\bigcup F^m_\gamma\cup\{\bigcup F^m_\gamma\}$. And $F^0_\gamma\subseteq
\bigcup F^0_\gamma\cup\{\bigcup F^0_\gamma\}$ therefore $\forall m(F^m_\gamma
\subseteq\bigcup F^m_\gamma\cup\{\bigcup F^m_\gamma\})$.\par
\begin{theorem}
If $\gamma$ is an ordinal then $\forall l(l\leq m\rightarrow F^m_\gamma
\subseteq F^l_\gamma)$.
\end{theorem}
{\it Proof.} $F^0_\gamma=\gamma$ and $\forall l(l\leq0\rightarrow l=0)$ thus
$\forall l(l\leq0\rightarrow F^0_\gamma\subseteq F^l_\gamma)$.\par
Assume if $x<m$ then $\forall l(l\leq x\rightarrow F^x_\gamma\subseteq
F^l_\gamma)$.\par
Then $F^{x+1}_\gamma\subseteq F^x_\gamma\wedge\forall l(l\leq x\rightarrow
F^x_\gamma\subseteq F^l_\gamma)$ thus $\forall l
(l\leq x\rightarrow F^{x+1}_\gamma\subseteq F^l_\gamma)$.\par
And $l=x+1\rightarrow F^{x+1}_\gamma\subseteq F^l_\gamma$ therefore $\forall l
(l\leq x+1\rightarrow F^{x+1}_\gamma\subseteq F^l_\gamma)$.\par
Apply the principle of mathematical induction to $0\leq x\wedge x<m$ and\par
conclude $\forall l(l\leq m\rightarrow F^m_\gamma\subseteq F^l_\gamma)$.
\begin{flushright}\boldmath$\Box$\end{flushright}
Then $F^m_\gamma\subseteq F^0_\gamma$ thus $\gamma\notin F^m_\gamma$. And
$\forall l(l\leq m\rightarrow\bigcup F^m_\gamma\subseteq\bigcup F^l_\gamma)$.\\
Apply the law of trichotomy and conclude $\bigcup F^m_\gamma\in\gamma
\leftrightarrow\bigcup F^m_\gamma\neq\gamma$.\\
Thus every first $\omega^m$-number is a first number, every first number is a
first $\omega^0$-number and $0$ is for every $m$ a first $\omega^m$-number. If
$\bigcup F^{m+1}_\alpha\neq\alpha$ then $\bigcup F^{m+1}_\alpha$ is the
greatest first $\omega^m$-number belonging to $\alpha$.\par
If $\bigcup F^m_\alpha=\alpha$ then $\forall l(l\leq m\rightarrow\alpha
\subseteq\bigcup F^l_\alpha\wedge\bigcup F^l_\alpha\subseteq\bigcup\alpha)$
thus $\bigcup\alpha=\alpha$ therefore $\forall l(l\leq m\rightarrow
\bigcup F^l_\alpha=\alpha)$. And $\forall u(u\in F^{m+1}_\alpha\leftrightarrow
u\in\alpha\wedge\bigcup F^m_u=u)$.\par
If $\bigcup F^m_\gamma\in F^m_\gamma$ then $\forall l(\bigcup F^{m+l}_\gamma
\neq\gamma)$ otherwise $\bigcup F^m_\gamma=\gamma$ contradictory to $\gamma
\notin F^m_\gamma$. Thus $\bigcup F^m_\gamma\in F^m_\gamma\rightarrow\forall l
(\bigcup F^{m+l}_\gamma\in\gamma)$.
\begin{axiom}
$\forall\alpha\forall m(\bigcup F^m_\alpha=\alpha\rightarrow\exists\gamma
\forall u(u\in\gamma\leftrightarrow\exists n(u\in\alpha+\omega^m\cdot n)))$.
\hfill\boldmath$\Box$\end{axiom}
$\gamma$ is denoted by $\alpha+\omega^{m+1}$. Thus $0+\omega^{m+1}$ is a set.\\
If $\bigcup\alpha=\alpha\wedge\beta\notin\alpha\wedge\beta\in\alpha+\omega^1$
then $\beta$ is an ordinal with first number $\alpha$. Thus every natural
number is a member of $0+\omega^1$ and every member of $0+\omega^1$ is a
natural number. This is the `axiom of infinity'.
\begin{definition}
If $\alpha$ is a first $\omega^m$-number then\par
$\alpha+\omega^{m+1}\cdot0=\alpha$,\par
$\alpha+\omega^{m+1}\cdot(n+1)=(\alpha+\omega^{m+1}\cdot n)+\omega^{m+1}$.
\hfill\boldmath$\Box$\end{definition}
If $\alpha$ is a first $\omega^{m+1}$-number then $\alpha$ is a first
$\omega^m$-number therefore\\
$\alpha+\omega^{m+1}=\alpha+\omega^{m+1}\cdot1$. And $\alpha+\omega^0=
\alpha+\omega^0\cdot1$ thus\par
$\forall m(\bigcup F^m_\alpha=\alpha\rightarrow\alpha+\omega^m=\alpha+\omega^m
\cdot1)$.
\begin{theorem}
$\bigcup F^m_\alpha=\alpha\rightarrow\alpha\in\alpha+\omega^m$.
\end{theorem}
{\it Proof.} If $\bigcup F^0_\alpha=\alpha$ then $\alpha+\omega^0=\alpha\cup
\{\alpha\}$ therefore $\alpha\in\alpha+\omega^0$.\par
Assume if $x<m$ then $\bigcup F^x_\alpha=\alpha\rightarrow\alpha\in\alpha+
\omega^x$.\par
If $\bigcup F^{x+1}_\alpha=\alpha$ then $\bigcup F^x_\alpha=\alpha$ thus
$\alpha\in\alpha+\omega^x$.\par
Then $\alpha+\omega^x=\alpha+\omega^x\cdot1$ thus $\alpha\in\alpha+\omega^x
\cdot1$ therefore $\alpha\in\alpha+\omega^{x+1}$.\par
Apply the principle of mathematical induction to $0\leq x\wedge x<m$ and\par
conclude $\bigcup F^m_\alpha=\alpha\rightarrow\alpha\in\alpha+\omega^m$.
\begin{flushright}\boldmath$\Box$\end{flushright}
\begin{theorem}
$\bigcup F^m_\alpha=\alpha\rightarrow\forall n(\bigcup F^m_{\alpha+\omega^{m+1}
\cdot n}=\alpha+\omega^{m+1}\cdot n)$.
\end{theorem}
{\it Proof.} If $\bigcup F^0_\alpha=\alpha$ then $\alpha+\omega^1$ exists
and\par
$\forall u(u\in\alpha+\omega^1\leftrightarrow\exists l(u\in\alpha+\omega^0\cdot
l))$. Then $\alpha+\omega^0\cdot l$ is an ordinal with\par
first number $\alpha$ thus a transitive set with transitive members.\par
Then $\bigcup(\alpha+\omega^1)\subseteq\alpha+\omega^1$ and
$\alpha+\omega^1\subseteq\bigcup(\alpha+\omega^1)$ thus $\alpha+\omega^1$ is a
first\par
$\omega^0$-number. Therefore $\bigcup F^0_\alpha=\alpha\rightarrow
\bigcup F^0_{\alpha+\omega^1}=\alpha+\omega^1$.\\
Assume if $y<l$ then $\alpha+\omega^1\cdot y$ is a first $\omega^0$-number.\par
Then $\alpha+\omega^1\cdot(y+1)=(\alpha+\omega^1\cdot y)+\omega^1$.\par
According to the induction assumption $\alpha+\omega^1\cdot y$ is a first
$\omega^0$-number.\par
Therefore $(\alpha+\omega^1\cdot y)+\omega^1$ is a first $\omega^0$-number thus
$\alpha+\omega^1\cdot(y+1)$ is\par
a first $\omega^0$-number.\\
Apply the principle of mathematical induction to $0\leq y\wedge y<l$ and\par
conclude $\bigcup F^0_\alpha=\alpha\rightarrow\forall l(\bigcup F^0_{\alpha+
\omega^{0+1}\cdot l}=\alpha+\omega^{0+1}\cdot l)$.\\
Assume if $x<m$ then $\bigcup F^x_\alpha=\alpha\rightarrow\forall n
(\bigcup F^x_{\alpha+\omega^{x+1}\cdot n}=\alpha+\omega^{x+1}\cdot n)$.\par
If $\alpha$ is a first $\omega^{x+1}$-number then apply the axiom of infinity
thus\par
$\exists\alpha+\omega^{(x+1)+1}\forall u(u\in\alpha+\omega^{(x+1)+1}
\leftrightarrow\exists l(u\in\alpha+\omega^{x+1}\cdot l))$.\par
$\alpha$ is a first $\omega^{x+1}$-number thus $\alpha$ is a first
$\omega^x$-number. Then according\par
to the induction assumption $\alpha+\omega^{x+1}\cdot l$ is a first
$\omega^x$-number.\par
If $v=\alpha+\omega^{x+1}\cdot l$ then $\bigcup F^x_v=v\wedge v\in v+
\omega^{x+1}$ therefore\par
$v\in\alpha+\omega^{x+1}\cdot(l+1)$ thus $v\in\alpha+\omega^{(x+1)+1}$.
Then\par
$\forall u(\exists v(u\in v\wedge v\in\alpha+\omega^{(x+1)+1}\wedge
\bigcup F^x_v=v)\rightarrow u\in\bigcup F^{x+1}_{\alpha+\omega^{(x+1)+1}})$.
\par
Then $\alpha+\omega^{(x+1)+1}\subseteq\bigcup F^{x+1}_{\alpha+
\omega^{(x+1)+1}}$. And $\bigcup F^{x+1}_{\alpha+\omega^{(x+1)+1}}\subseteq
\alpha+\omega^{(x+1)+1}$.\par
Therefore $\bigcup F^{x+1}_{\alpha+\omega^{(x+1)+1}}=\alpha+
\omega^{(x+1)+1}$.\par
Thus $\bigcup F^{x+1}_\alpha=\alpha\rightarrow\bigcup F^{x+1}_{\alpha+
\omega^{(x+1)+1}}=\alpha+\omega^{(x+1)+1}$.\\
Assume if $y<l$ then $\alpha+\omega^{(x+1)+1}\cdot y$ is a first
$\omega^{x+1}$-number.\par
Then $\alpha+\omega^{(x+1)+1}\cdot(y+1)=(\alpha+\omega^{(x+1)+1}\cdot y)+
\omega^{(x+1)+1}$.\par
According to the induction assumption $\alpha+\omega^{(x+1)+1}\cdot y$ is a
first\par
$\omega^{x+1}$-number. Therefore $(\alpha+\omega^{(x+1)+1}\cdot y)+
\omega^{(x+1)+1}$ is a first $\omega^{x+1}$-number\par
thus $\alpha+\omega^{(x+1)+1}\cdot(y+1)$ is a first $\omega^{x+1}$-number.\\
Apply the principle of mathematical induction to $0\leq y\wedge y<l$ and\par
conclude $\bigcup F^{x+1}_\alpha=\alpha\rightarrow\forall l
(\bigcup F^{x+1}_{\alpha+\omega^{(x+1)+1}\cdot l}=\alpha+\omega^{(x+1)+1}\cdot
l)$.\\
Apply the principle of mathematical induction to $0\leq x\wedge x<m$ and\par
conclude $\bigcup F^m_\alpha=\alpha\rightarrow\forall n
(\bigcup F^m_{\alpha+\omega^{m+1}\cdot n}=\alpha+\omega^{m+1}\cdot n)$.
\begin{flushright}\boldmath$\Box$\end{flushright}
With the principle of mathematical induction one can prove\par
$\forall n\forall l(\bigcup F^m_\alpha=\alpha\rightarrow(\alpha+\omega^{m+1}
\cdot n)+\omega^{m+1}\cdot l=\alpha+\omega^{m+1}\cdot(n+l))$.
\begin{theorem}
$\bigcup F^m_\alpha=\alpha\rightarrow F^{m+1}_{\alpha+\omega^{m+1}}=
F^{m+1}_\alpha\cup\{\alpha\}$.
\end{theorem}
{\it Proof.} If $\alpha$ is a first $\omega^0$-number then $\bigcup\alpha=
\alpha$ and\par
$\forall u(u\in F^1_{\alpha+\omega^1}\leftrightarrow u\in\alpha+\omega^1\wedge
\bigcup u=u)$.\par
And $\forall u(u\in\alpha+\omega^1\leftrightarrow\exists n(u\in\alpha+\omega^0
\cdot n))$.\par
Then $\alpha+\omega^0\cdot n$ is an ordinal with first number $\alpha$
therefore\par
$\forall u(u\in\alpha+\omega^0\cdot n\cup\{\alpha+\omega^0\cdot n\}\rightarrow
u\in\alpha\cup\{\alpha\}\vee\bigcup u\in u)$.\par
Thus $\forall u(u\in F^1_{\alpha+\omega^1}\rightarrow(u\in\alpha\cup\{\alpha\}
\vee\bigcup u\in u)\wedge\bigcup u=u)$.\par
Then $F^1_{\alpha+\omega^1}\subseteq F^1_\alpha\cup\{\alpha\}$. And
$F^1_\alpha\cup\{\alpha\}\subseteq F^1_{\alpha+\omega^1}$ thus
$F^1_{\alpha+\omega^1}=F^1_\alpha\cup\{\alpha\}$.\par
Therefore $\bigcup F^0_\alpha=\alpha\rightarrow F^{0+1}_{\alpha+\omega^{0+1}}=
F^{0+1}_\alpha\cup\{\alpha\}$. Then $\bigcup F^{0+1}_{\alpha+\omega^{0+1}}=
\alpha$.\par
Thus $\alpha$ is the greatest first $\omega^0$-number belonging to $\alpha+
\omega^{0+1}$.\\
Assume if $x<m$ then $\bigcup F^x_\alpha=\alpha\rightarrow F^{x+1}_{\alpha+
\omega^{x+1}}=F^{x+1}_\alpha\cup\{\alpha\}$.\par
Then $\alpha$ is the greatest first $\omega^x$-number belonging to $\alpha+
\omega^{x+1}$.\par
If $\alpha$ is a first $\omega^{x+1}$-number then apply the axiom of infinity
thus\par
$\exists\alpha+\omega^{(x+1)+1}\forall u(u\in\alpha+\omega^{(x+1)+1}
\leftrightarrow\exists l(u\in\alpha+\omega^{x+1}\cdot l))$.\par
Then $\forall u(u\in F^{(x+1)+1}_{\alpha+\omega^{(x+1)+1}}\leftrightarrow
u\in\alpha+\omega^{(x+1)+1}\wedge\bigcup F^{x+1}_u=u)$.\par
And $\forall u(u\in\alpha+\omega^{(x+1)+1}\leftrightarrow\exists l(u\in\alpha
\cup\{\alpha\}\vee u\in\alpha+\omega^{x+1}\cdot(l+1)))$.\par
Then $\alpha$ is a first $\omega^{x+1}$-number thus $\alpha$ is a first
$\omega^x$-number and $\alpha+\omega^{x+1}\cdot l$\par
is a first $\omega^x$-number. Then according to the induction assumption\par
$F^{x+1}_{(\alpha+\omega^{x+1}\cdot l)+\omega^{x+1}}=F^{x+1}_{\alpha+
\omega^{x+1}\cdot l}\cup\{\alpha+\omega^{x+1}\cdot l\}$. Thus $\alpha+
\omega^{x+1}\cdot l$ is the\par
greatest first $\omega^x$-number belonging to $\alpha+\omega^{x+1}\cdot
(l+1))$.\par
Then $\forall u(u\in\alpha+\omega^{(x+1)+1}\wedge\bigcup F^{x+1}_u=u\rightarrow
u\in\alpha\cup\{\alpha\}\wedge\bigcup F^{x+1}_u=u)$.\par
Thus $F^{(x+1)+1}_{\alpha+\omega^{(x+1)+1}}\subseteq F^{(x+1)+1}_\alpha\cup
\{\alpha\}$. And $F^{(x+1)+1}_\alpha\cup\{\alpha\}\subseteq
F^{(x+1)+1}_{\alpha+\omega^{(x+1)+1}}$.\par
Therefore $F^{(x+1)+1}_{\alpha+\omega^{(x+1)+1}}=F^{(x+1)+1}_\alpha\cup
\{\alpha\}$.\\
Apply the principle of mathematical induction to $0\leq x\wedge x<m$ and\par
conclude $\bigcup F^m_\alpha=\alpha\rightarrow F^{m+1}_{\alpha+\omega^{m+1}}=
F^{m+1}_\alpha\cup\{\alpha\}$. Then $\bigcup F^{m+1}_{\alpha+\omega^{m+1}}=
\alpha$.\par
Therefore $\alpha$ is the greatest first $\omega^m$-number belonging to
$\alpha+\omega^{m+1}$.
\begin{flushright}\boldmath$\Box$\end{flushright}
\begin{theorem}
$\bigcup F^m_\gamma=\gamma\wedge\bigcup F^{m+1}_\gamma\in\gamma\rightarrow
\gamma=\bigcup F^{m+1}_\gamma+\omega^{m+1}$.
\end{theorem}
{\it Proof.} $\bigcup F^{m+1}_\gamma\in\gamma$ thus $\bigcup F^{m+1}_\gamma$ is
the greatest first $\omega^m$-number belonging to\par
$\gamma$. And $\bigcup F^{m+1}_\gamma+\omega^{m+1}$ is a first
$\omega^m$-number. Thus $\bigcup F^{m+1}_\gamma$ is the greatest\par
first $\omega^m$-number belonging to $\bigcup F^{m+1}_\gamma+\omega^{m+1}$.\par
Then $\bigcup F^{m+1}_\gamma+\omega^{m+1}\notin\gamma\wedge\gamma\notin
\bigcup F^{m+1}_\gamma+\omega^{m+1}$.\par
Apply the law of trichotomy and conclude $\gamma=\bigcup F^{m+1}_\gamma+
\omega^{m+1}$.
\begin{flushright}\boldmath$\Box$\end{flushright}
\begin{theorem}
$\bigcup F^m_\gamma=\gamma\rightarrow\exists\alpha\exists n
(\bigcup F^{m+1}_\alpha=\alpha\wedge\gamma=\alpha+\omega^{m+1}\cdot n)$.
\end{theorem}
{\it Proof.} $\bigcup F^{m+1}_\gamma=\gamma\rightarrow\bigcup F^m_\gamma=
\gamma\wedge\gamma=\gamma+\omega^{m+1}\cdot0$ thus $\alpha=\gamma$ and
$n=0$.\par
If $\bigcup F^{m+1}_\gamma\neq\gamma$ then $\bigcup F^{m+1}_\gamma\in
\gamma$ thus $\gamma=\bigcup F^{m+1}_\gamma+\omega^{m+1}\cdot1$.\par
Use separation to define set $v$ of first $\omega^m$-numbers by\par
$\forall u(u\in v\leftrightarrow u\in F^{m+1}_\gamma\wedge u\notin u\wedge
\exists n(\gamma=u+\omega^{m+1}\cdot n))$. Then $\bigcup F^{m+1}_\gamma\in
v$.\par
Apply regularity to $v$ to find first $\omega^m$-number $\alpha\in v$ such
that\par
$\alpha\in F^{m+1}_\gamma\wedge\exists n(\gamma=\alpha+\omega^{m+1}\cdot n)
\wedge\forall u(u\in\alpha\wedge u\notin u\rightarrow u\notin v)$.\par
Then $\bigcup F^{m+1}_\alpha\in\alpha\cup\{\alpha\}$.\par
If $\bigcup F^{m+1}_\alpha\in\alpha$ then $\bigcup F^m_\alpha=\alpha\wedge
\bigcup F^{m+1}_\alpha\in\alpha\rightarrow\alpha=\bigcup F^{m+1}_\alpha+
\omega^{m+1}\cdot1$.\par
Thus $\bigcup F^{m+1}_\alpha$ is the greatest first $\omega^m$-number belonging
to $\alpha$ and\par
$\gamma=\bigcup F^{m+1}_\alpha+\omega^{m+1}\cdot(1+n)$.\par
Then $\bigcup F^{m+1}_\alpha\in\alpha\wedge\alpha\in F^{m+1}_\gamma$ thus
$\alpha\in\gamma$ therefore $\bigcup F^{m+1}_\alpha\in F^{m+1}_\gamma$.\par
Then $\bigcup F^{m+1}_\alpha\in v$ contradictory to $\forall u(u\in\alpha\wedge
u\notin u\rightarrow u\notin v)$.\par
Therefore $\bigcup F^m_\gamma=\gamma\rightarrow\exists\alpha\exists n(\bigcup
F^{m+1}_\alpha=\alpha\wedge\gamma=\alpha+\omega^{m+1}\cdot n)$.
\begin{flushright}\boldmath$\Box$\end{flushright}
If $\bigcup F^m_\alpha=\alpha$ then $\forall n(\bigcup F^m_{\alpha+\omega^{m+1}
\cdot n}=\alpha+\omega^{m+1}\cdot n)$ thus $\alpha+\omega^{m+1}\cdot n$ is the
greatest first $\omega^m$-number belonging to $\alpha+\omega^{m+1}\cdot(n+1)$.
If $\bigcup F^{m+1}_\alpha=\alpha$ then for every $n\neq0$ is $\alpha$ the
greatest first $\omega^{m+1}$-number belonging to $\alpha+\omega^{m+1}
\cdot n$.\par
Every first number is a first $\omega^0$-number thus every ordinal is an
ordinal with a first $\omega^0$-number. We define `$\omega^{m+1}$-numbers with
first number $\alpha$'.
\begin{definition}
Ordinal $\gamma$ is an `$\omega^{m+1}$-number with first number $\alpha$'
if\par
$\bigcup F^m_\gamma=\gamma\wedge\bigcup F^{m+1}_\alpha=\alpha\wedge$\par
$\mbox{ }\mbox{ }\mbox{ }\gamma\notin\alpha\wedge\forall u(u\in\gamma\cup
\{\gamma\}\rightarrow u\in\alpha\cup\{\alpha\}\vee\bigcup F^{m+1}_u\in u)$.\par
The successor of $\omega^{m+1}$-number $\gamma$ with first number $\alpha$ is
$\gamma+\omega^{m+1}$.
\hfill\boldmath$\Box$\end{definition}
If $\gamma\neq\alpha$ then $\gamma=\bigcup F^{m+1}_\gamma+\omega^{m+1}$ thus
$\exists\beta\exists n(\bigcup F^{m+1}_\beta=\beta\wedge\gamma=\beta+\omega^
{m+1}\cdot n)$. Then $\gamma\neq\beta$ thus $n\neq0$ and $\beta$ is the
greatest first $\omega^{m+1}$-number belonging to $\gamma$. Thus $\alpha\in
\beta\cup\{\beta\}$. And $\beta\in\gamma\wedge\bigcup F^{m+1}_\beta=\beta
\rightarrow\beta\in\alpha\cup\{\alpha\}$ thus $\beta=\alpha$. If $\gamma=
\alpha$ then $\gamma=\alpha+\omega^{m+1}\cdot0$. Thus $\exists n(\gamma=\alpha+
\omega^{m+1}\cdot n)$.\par
Ordinals with a first $\omega^0$-number satisfy the Peano axioms. By the same
kind of reasoning used in section 5 about natural numbers one can prove if
$\gamma$ is an $\omega^{m+1}$-number with first number $\alpha$ then
\begin{itemize}
\item there is first $\omega^{m+1}$-number $\alpha$,
\item the successor $\gamma+\omega^{m+1}$ of $\gamma$ is an
$\omega^{m+1}$-number with first number $\alpha$,
\item $\alpha$ is not the successor of an $\omega^{m+1}$-number with first
number $\alpha$,
\item if $\beta\notin\alpha\wedge\beta\in\gamma\wedge\bigcup F^m_\beta=\beta$
then $\beta$ is an $\omega^{m+1}$-number with first number $\alpha$,
\item if $\beta$ is an $\omega^{m+1}$-number with first number $\alpha$ and
$\beta+\omega^{m+1}=\gamma+\omega^{m+1}$ then $\beta\in\gamma\cup\{\gamma\}
\wedge\gamma\in\beta\cup\{\beta\}$ thus $\beta=\gamma$,
\item if $\Phi(\beta)$ is a well formed formula with first $\omega^m$-number
$\beta$ free in $\Phi$ then\\
$\Phi(\alpha)\wedge\forall\beta(\beta\notin\alpha\wedge\beta\in\gamma\wedge
\bigcup F^m_\beta=\beta\wedge\Phi(\beta)\rightarrow\Phi(\beta+\omega^{m+1}))
\rightarrow\Phi(\gamma)$.
\end{itemize}
Thus $\omega^m$-numbers with first number $\alpha$ satisfy the Peano axioms.
Therefore\\
\mbox{ }\par
if $\beta$ is an ordinal then $\exists\alpha\exists n(\bigcup\alpha=\alpha
\wedge\beta=\alpha+\omega^0\cdot n)$,\par
if $\gamma$ is a first $\omega^m$-number then $\exists\alpha\exists n(\bigcup
F^{m+1}_\alpha=\alpha\wedge\gamma=\alpha+\omega^{m+1}\cdot n)$.
\section{Counting down ordinals}
Ordinal $\alpha$ is a `first $\omega^\omega$-number' if $\forall m
(\bigcup F^m_\alpha=\alpha)$. $0$ is a first $\omega^\omega$-number.\\
If $\alpha$ is a first $\omega^\omega$-number then for every $m$ and for every $n$ $(\alpha+\omega^{m+1}\cdot
n)$ is an $\omega^{m+1}$-number with first number $\alpha$.
\begin{theorem}
If $\gamma$ is an ordinal then $\exists\alpha(\forall m(\bigcup F^m_\alpha=
\alpha)\wedge\exists l(\bigcup F^{l+1}_\gamma=\alpha))$.
\end{theorem}
{\it Proof.} If ordinal $\gamma$ is a first $\omega^\omega$-number then
$\alpha=\gamma$ thus $\bigcup F^{0+1}_\gamma=\alpha$.\\
If $\gamma$ is not a first $\omega^\omega$-number then $\exists n
(\bigcup F^n_\gamma\in\gamma)$ thus\par
$\bigcup F^0_\gamma\in\gamma$ or $\exists m(m<n\wedge\bigcup F^m_\gamma=\gamma
\wedge\bigcup F^{m+1}_\gamma\in\gamma)$.\par
If $\bigcup F^0_\gamma\in\gamma$ then $\gamma$ is an ordinal with a first
number.\par
Therefore $\bigcup F^0_\gamma\in\gamma\wedge\exists\beta\exists n
(\bigcup F^0_\beta=\beta\wedge n\neq0\wedge\gamma=\beta+\omega^0\cdot n\wedge
\bigcup F^{0+1}_\gamma=\beta)$.\par
If $\exists m(\bigcup F^m_\gamma=\gamma\wedge\bigcup F^{m+1}_\gamma\in\gamma)$
then\par
$\exists\beta\exists n(\bigcup F^{m+1}_\beta=\beta\wedge n\neq0\wedge
\gamma=\beta+\omega^{m+1}\cdot n\wedge\bigcup F^{(m+1)+1}_\gamma=\beta)$.\\
Thus in both cases $\exists\beta\exists n$ such that we can `countdown $\gamma$
to $\beta$ in $n$ steps'.\par
If $\forall l(\bigcup F^l_\beta=\beta)$ then $\beta$ is a first
$\omega^\omega$-number and $\beta\in\gamma$.\par
If $\exists m(\bigcup F^m_\beta=\beta\wedge\bigcup F^{m+1}_\beta\in\beta)$ then
$\beta=\bigcup F^{m+1}_\beta+\omega^{m+1}$. Therefore\par
$\exists\beta_1(\bigcup F^{m+1}_{\beta_1}=\beta_1\wedge\exists n_1(n_1\neq0
\wedge\beta=\beta_1+\omega^{m+1}\cdot n_1))$. Thus\par
$\mbox{ }\gamma=(\beta_1+\omega^{m+1}\cdot n_1)+\omega^0\cdot n$. Then
$\beta_1$ is the greatest first $\omega^{m+1}$-number\par
belonging to $\beta$ and $\beta\in\gamma$ therefore $\bigcup F^{(m+1)+1}_\gamma
=\beta_1$.\\
Thus we can countdown $\gamma$ to $\beta_1$ in $n+n_1$ steps.\par
If $\forall l(\bigcup F^l_{\beta_1}=\beta_1)$ then $\beta_1$ is a first
$\omega^\omega$-number and $\beta_1\in\gamma$.\par
If $\exists m_1(\bigcup F^{m_1}_{\beta_1}=\beta_1\wedge
\bigcup F^{m_1+1}_{\beta_1}\in\beta_1)$ then $\beta_1=
\bigcup F^{m_1+1}_{\beta_1}+\omega^{m_1+1}$. Therefore\par
$\exists\beta_2(\bigcup F^{m_1+1}_{\beta_2}=\beta_2\wedge\exists n_2(n_2\neq0
\wedge\beta_1=\beta_2+\omega^{m_1+1}\cdot n_2))$. Thus\par
$\mbox{ }\gamma=((\beta_2+\omega^{m_1+1}\cdot n_2)+\omega^{m+1}\cdot n_1)+
\omega^0\cdot n$. Then $\beta_2$ is the greatest first\par
$\omega^{m_1+1}$-number belonging to $\beta_1$ and $\beta_1\in\gamma$ therefore
$\bigcup F^{(m_1+1)+1}_\gamma=\beta_2$.\\
Thus we can countdown $\gamma$ to $\beta_2$ in $n+n_1+n_2$ steps.\par
If $\forall l(\bigcup F^l_{\beta_2}=\beta_2)$ then $\beta_2$ is a first
$\omega^\omega$-number and $\beta_2\in\gamma$.\\
Thus there is a rule to construct a descending sequence\\
$\mbox{ }\mbox{ }\mbox{ }\mbox{ }\ldots,\mbox{ }\beta_2\in\beta_1,\mbox{ }
\beta_1\in\beta,\mbox{ }\beta\in\gamma$.\\
If $\beta$ is a first $\omega^\omega$-number then we define $\alpha=\beta$
thus\par
$\exists\alpha(\forall m(\bigcup F^m_\alpha=\alpha)\wedge\exists l
(\bigcup F^{l+1}_\gamma=\alpha))$.\\
If $\beta$ is not a first $\omega^\omega$-number then we define $\beta_0=
\beta$.\par
Assume $\forall n\exists\beta_{n+1}(\beta_{n+1}\in\beta_n)$. Thus $\beta_{n+1}
\in\beta\wedge\beta_{n+1}\in\gamma$.\par
Use separation to define set $s$ by\par
$\forall u(u\in s\leftrightarrow u\in\gamma\wedge\exists n(u=\beta_n))$ thus
$s\subseteq\gamma$ and set $\beta_0\in s$. Then\par
$\forall v(v\in s\rightarrow\exists n(v=\beta_n\wedge\beta_{n+1}\in v\wedge
\beta_{n+1}\in s))$ contradictory to regularity\par
$\exists v(v\in s\wedge v\notin v\wedge\forall u(u\in v\wedge u\notin u
\rightarrow u\notin s))$.\\
Therefore $\exists\alpha(\forall m(\bigcup F^m_\alpha=\alpha)\wedge\exists l
(\bigcup F^{l+1}_\gamma=\alpha))$.
\begin{flushright}\boldmath$\Box$\end{flushright}
If $\gamma\neq\alpha$ then $\alpha$ is the greatest first $\omega^
\omega$-number belonging to $\gamma$. Then counting down $\gamma$ terminates at
$\alpha$ in a finite number of steps. In this way, with the Peano axioms, if
$\alpha\neq0$ then $\alpha$ is an impassable barrier for counting down $\gamma$
to $0$ in a finite number of steps.\par
If $\gamma=0$ or $0$ is the greatest first $\omega^\omega$-number belonging to
$\gamma$ then counting down $\gamma$ terminates at $0$ in a finite number of
steps. And if $\gamma\neq0$ then $\exists l(\bigcup F^{l+1}_\gamma=0\wedge
F^{l+1}_\gamma=1)$.
\section{References}
\begin{description}
\item P.R. Halmos. \em Naive Set Theory \em\\
Springer-Verlag New York Inc., 1974
\item Thomas Jech. \em Set Theory\em\\
The Third Millennium Edition, revised and expanded, May 2002\\
Springer Monographs in Mathematics
\item Elliott Mendelson. \em Introduction to Mathematical Logic\em\\
Sixth Edition, 2015
\end{description}
\end{document}